# Event-triggered Communication in Wide-area Damping Control: A Limited Output Feedback Based Approach

Mahendra Bhadu, *Student Member, IEEE*, Niladri Sekhar Tripathy, I. N. Kar, *Senior Member, IEEE,* and Nilanjan Senroy, *Member, IEEE*

***Abstract*—A conceptual design methodology is proposed for event-triggered based power system wide area damping controller. The event-triggering mechanism is adopted to reduce the communication burden between origin of the remote signal and the wide area damping controller (WADC) location. The remote signal is transmitted to the WADC only when an event-triggering condition based on a predefined system output, is satisfied. The triggering condition is derived from a stability criterion, and is monitored continuously by a separate event-monitoring unit located at the origin of the remote signal. The stability of the resulting closed loop system is guaranteed via the input-to-state stability (ISS) technique. The proposed event triggered WADC (ET-WADC) is implemented on two typical test power systems - two area four machine and IEEE 39 bus 10 machine. The validation of proposed mechanism is carried out through non-linear simulation studies on MATLAB/Simulink platform. The numerical results show the efficacy of the controller in managing the communication channel usage without compromising the stated system stability objectives.**

***Index Terms*— Aperiodic control, event-triggered control, inter-area oscillations, input-to-state stability, output-feedback control, power system control, wide area damping controller.**

## I. Introduction

MODERN power system networks are complex and widely distributed in nature. Any kind of disturbance ignites electromechanical oscillations, which if not properly damped can pose major challenges in system operation. Historically, the advent of power system stabilizers (PSS) alleviated the power oscillation problems in networks where large amounts of power were transmitted over weak links. In more recent times, spatially distributed electromechanical oscillations in the frequency range of 0.2-0.8 Hz have proved to be difficult to damp out using the conventional PSS design. With the advent of geographically distributed synchronized phasor measurement units (PMUs), it has been proposed to use remote signals to drive damping controllers so as to lessen the impact of such inter-area oscillations. The wide area damping controller (WADC) is an imperative tool in modern power system for small signal and transient stability enhancement [1]-[4]. Worldwide, these WADC are still at testing stage primarily because it is not always possible to ensure the availability of the remote signal for the controller, due to various communication network issues.

The authors are with the Department of Electrical Engineering, Indian Institute of Technology (IIT) Delhi, New Delhi 110016, India.(e-mails: mbhadu@gmail.com, niladri.tripathy@ee.iitd.ac.in, ink@ee.iitd.ac.in and nsenroy@ee.iitd.ac.in).

The design and application of WADC in the presence of imperfect communication medium has been studied in recent years [2], [5]-[9]. Most proposed control schemes rely on the continuous availability of the remote signal for their success. Periodic updating of the measurement at the controller location, is based on a pessimistic control approach, where worst case scenarios are considered while designing the measurement signal sampling strategy. Such scenarios are rare and consequently the communication medium may be unnecessarily congested even when the power system is stable and requires little attention. In case of shared communication or limited bandwidth availability of the medium, the problem becomes acute. Hence, it is advantageous to reduce the communication burden between the sensors and controller location in the given network controlled system (NCS) [10]-[12]. This is a motivation to adopt an aperiodic transmission technique over the continuous one.

The benefits of aperiodic sampling are reported in [13]-[15]. Recently, [14]-[16] have proposed event based aperiodic control technique, where control signal is computed and transmitted at the system end only when a predefined event condition is violated. Other researchers have use this event-triggering technique to solve control, estimation and communication problems [17]-[19]. In [20] Zhang et al proposed an event triggering based wide area damping control technique in power systems. The event-triggering condition depends on system states, which may not necessarily be observable at all times. The present problem is formulated differently i.e. an output dependent event triggering condition is proposed in this paper, that depends on a generator speed information, which is realizable in practical scenario. The stability condition, derived for an output feedback system, is a more challenging control problem as compared to state feedback. The proposed ET-WADC is designed to lower the communication burden while maintaining a minimum damping effect on inter-area electromechanical oscillations.

The ET-WADC updates its remote signal input only when an event-triggering condition is satisfied. The tuning parameter of the event triggering mechanism is adjusted to obtain the satisfactory damping effect on one hand and reduced communication on the other. The stability of entire closed loop system is ensured through (input-to-state stability) ISS technique. The generic form of output feedback based event triggered control mechanism is depicted along-with one of the test power system, in Fig. 1.

The selection of WADC locations and input signals is carried out on the basis of the established modal residue method [21]-[23]. Time delays are typically associated with PMU to WADC data transmission, due to measurement and distances involved in communication medium as well as network congestion. The data processing and alignment time taken by phasor data concentrators is typically 2 *ms* - 2 *s* [24]. The net data latency has an adverse effect on system stability and the damping effectiveness of the WADC [25]-[29]. For the purpose of controller design, this time delay is usually represented using the Pade approximation technique [25]. In this paper, a constant time delay is 0.1 *s* is selected. The equivalent third order Pade approximation for a delay of *T* second is as given in (1).

$$TF_{\text{Pade}} \approx \frac{-\frac{1}{120}s^3T^3 + \frac{1}{12}s^2T^2 - \frac{1}{2}sT + 1}{\frac{1}{120}s^3T^3 + \frac{1}{12}s^2T^2 + \frac{1}{2}sT + 1} \tag{1}$$

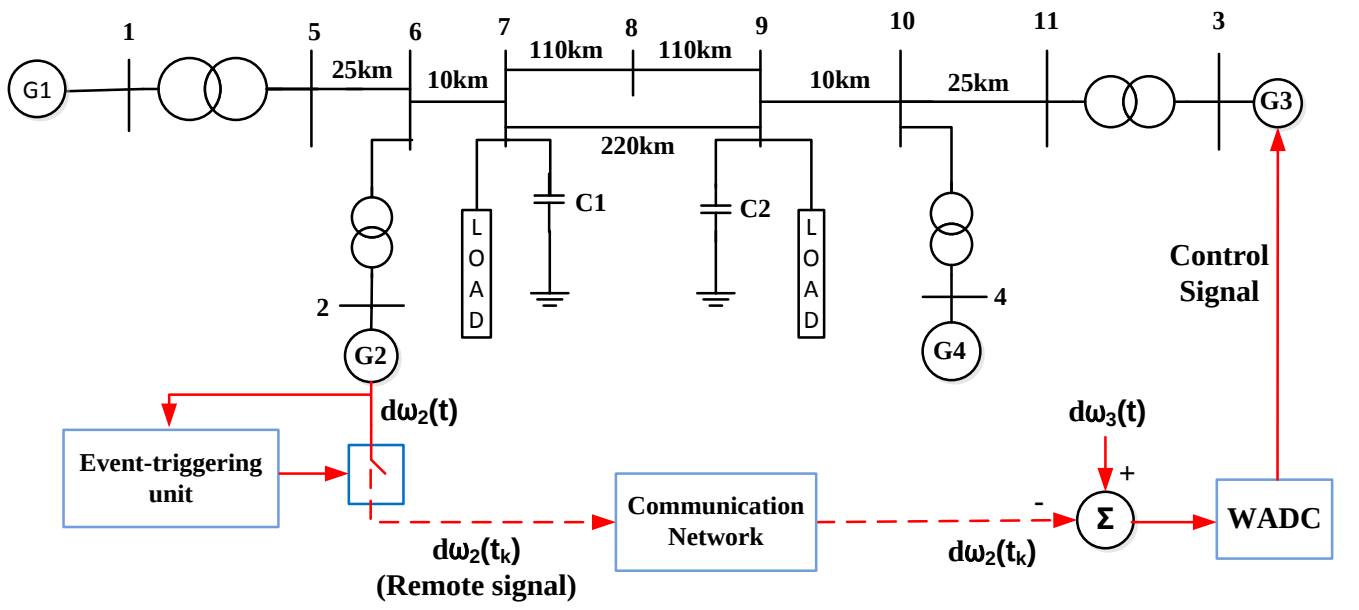

**Fig. 1.** Conceptual Block diagram of event-triggered based two-area four-machine power system. Here the bold and dotted line represent the continuous and aperiodic communication link respectively.

*Organization*

The rest of the paper is structured as follows. Section II presents modelling of event-triggered communication in wide area control . Case studies of multi-machine power systems are discussed in section IV, which also includes the methodology for selection of controller location and it's input signals. Time domain simulations results obtained using the event triggered controller with actual nonlinear power system models are presented in Section V. Section VI is the concluding section.

*A. Notation*

Throughout this paper the symbols $\mathbb{R}$,$\mathbb{R}^+$ and *I* are used to denote a set of real numbers, positive real numbers and non-negative integers respectively. The Euclidean norm of any vector $\mathbb{R}^n$ is represented by the symbol ||.||, where superscript *n* indicates the dimension of that vector. Symbol $\lambda_{min}(P)$ of any real matrix $P^{n\times m}$ with row and column dimension *n* and *m*, is used to define the minimum eigen value of matrix *P*. The notation $P > 0$, $P < 0$, $P^{-1}$ and $P^T$ represents the positive definiteness, negative definiteness, inverse and transpose of a given matrix *P*. The symbol '*inf*' and '$\wedge$' denote the infimum of any sequence and logical AND operation accordingly. Apart from the above defined symbols, the following definitions are also essential to understand theoretical concepts:

**Definition 1.** A continuous function $f:[0,a) \rightarrow [0\ \infty)$ is a class $K_\infty$ function if it holds the following conditions [30]:

1. $f(0) = 0$ *when* $a = 0$
2. $f(a) \rightarrow \infty$ *when* $a \rightarrow \infty$
3. $df/da > 0\ \forall\ a > 0$

**Definition 2.** A system

$$\dot{x} = f(x,e) \tag{2}$$

is said to input-to-state stable (ISS) [31],[32], with respect to its input $e$, if there exist a continuous time differentiable ISS Lyapunov function $V(x): \mathbb{R}^n \rightarrow \mathbb{R}$, which satisfy the following conditions:

$$\alpha_1(||x||) \leq V(x) \leq \alpha_2(||x||) \tag{3}$$

$$\dot{V}(x) \leq -\alpha_3(||x|| + \gamma(||e||)) \tag{4}$$

where $\alpha_1, \alpha_2, \alpha_3$ and $\gamma$ are class $K_\infty$ functions.

## II. MODELLING OF EVENT-TRIGGERED COMMUNICATION IN WIDE AREA CONTROL

This section introduces a novel control strategy for network controlled power systems with minimum communication network usage. The key idea of the present strategy is, the communication network will be used only when a predefined condition is violated. This predefined condition is called as an event-triggering condition, derived from the system stability criteria. *The event refers to the transmission of the remote signal to the controller location*. Let us consider a continuous-time linear two area power system model as:

$$\dot{x}_p(t) = A_p x_p(t) + B_p u_p(t) \tag{5}$$

$$y_1(t) = C_1 x_p(t),\ y_1 \in \mathbb{R} \tag{6}$$

$$y_2(t) = C_2 x_p(t),\ y_2 \in \mathbb{R} \tag{7}$$

here, $x_p \in \mathbb{R}^{n\times 1}$ and $u_p \in \mathbb{R}^{m\times 1}$ are system states and control input respectively. The dimension of system matrices are $A_p \in \mathbb{R}^{n\times n}$, $B_p \in \mathbb{R}^{n\times m}$, $C_1 \in \mathbb{R}^{m\times n}$ and $C_2 \in \mathbb{R}^{m\times n}$. The system have two outputs as $y_1 = d\omega_2$ and $y_2 = d\omega_3$. Here $d\omega_2$ and $d\omega_3$ denotes the speed deviation of the machine 2 and 3 respectively, with respect to machine 4. To control (5), a input $u_p(t)$ is generated through dynamic WADC as mentioned in (8) and (9) as:

$$\dot{x}_c(t) = A_c x_c(t) + B_c(y_2(t) - y_1(t)) \tag{8}$$

$$u_p(t) = C_c x_c(t) + D_c(C_2 - C_1)x_p(t) \tag{9}$$

here, $x_c$ represents the controller states and the difference between two system outputs $(y_2 - y_1)$ acts as a input to the controller. The system (5) and controller (8) can be written as a continuous-time augmented system as in (10).

$$\begin{bmatrix}\dot{x}_p \\ \dot{x}_c\end{bmatrix} = \begin{bmatrix} A_p + B_p D_c C & B_p C_c \\ B_c C & A_c \end{bmatrix} \begin{bmatrix} x_p \\ x_c \end{bmatrix} \tag{10}$$

where, matrix $C = (C_2 - C_1)$. According to Fig. 1, the plant output $y_1$ is not locally available at controller location and it is transmitted eventually via a communication link. Suppose $t_{k\epsilon I}$, represents the event-triggering sequences at which the output $y_1$ is transmitted to the controller through the communication link. Then the system and controller dynamics (5-9) using the eventual output $y_1(t_k)$ can be written as [14], [16]:

$$\dot{x}_p(t) = A_p x_p(t) + B_p u_p(t) \tag{11}$$

$$\dot{x}_c(t) = B_c C x_p(t) + A_c x_c(t) - B_c e_y(t) \tag{12}$$

$$u_p(t) = C_c x_c(t) + D_c C x_p(t) - D_c e_y(t) \tag{13}$$

The variable $e_y(t) \in \mathbb{R}$ is described as measurement error of output $y_1(t)$ and defined as:

$$e_y(t) = y_1(t_k) - y_1(t), \forall\, t \in [t_k, t_{k+1}) \tag{14}$$

Using (11)-(13), the event-triggered augmented system is:

$$\dot{x}(t) = Ax(t) + Be_y(t) \tag{15}$$

where $\dot{x}(t) = \begin{bmatrix} \dot{x}_p \\ \dot{x}_c \end{bmatrix}$, $A = \begin{bmatrix} A_p + B_p D_c C & B_p C_c \\ B_c C & A_c \end{bmatrix}$, $B = \begin{bmatrix} -B_c \\ -B_p D_c \end{bmatrix}$.

### *A. Problem Description and Statement*

The primary aim of this paper is to design an output feedback based event-triggering law for WADC. To achieve the control input, the outputs $y_1$ and $y_2$ are used as input signals of the WADC. The WADC is co-located with the generator 3, therefore WADC can access the output signal $y_2$ continuously but the signal $y_1$ (speed deviation, $d\omega_2$) needs to be transmitted via communication medium as shown in Fig. 1. The continuous transmission of $y_1$ to the WADC means the maximum usage of communication bandwidth. Generally, communication channels are shared in nature, therefore continuous availability of channel is sometime quite impossible or economically costly. To resolve this issue the present approach proposes a novel aperiodic-control technique in which continuous transmission of output signal $y_1$ is not required. The signal $y_1$ is transmitted to the WADC from the generator 2, only when the $||y_1||$ exceeds a threshold, to make the closed loop system stable. This predefined output dependent threshold condition is called as event-triggering condition.

**Problem Statement:** To design an event-triggering condition and the controller gain $C_c$ and $D_c$ of (13), which will ensure the ISS of (11), with respect to measurement error $e_y(t)$.

The above mentioned problem has been solved using following steps, (i) controller design, (ii) event-triggering condition design.

### *B. Controller Design*

In order to better illustrate the event triggering concept, a classical lead-lag type controller is used [33]. This may be extended to any robust or adaptive controller. The transfer function of the WADC is given by:

$$TF_{WADC}(s) = K\left(\frac{sT_w}{sT_w + 1}\right)\left(\frac{1 + s\tau_1}{1 + s\tau_2}\right) \tag{16}$$

where $T_w$ is the washout time constant, typically 10 *s*. The $\tau_1$ and $\tau_2$ are time constants of the compensating network. The gain *K* was chosen to provide sufficient damping for inter area mode of oscillations.

The above designed controller (12) is used to prove the stability of (11) through ISS theory [31], [32], [41]. A brief idea on ISS is stated in the form of definition in previous section. The design procedure of event-triggering condition is discussed elaborately in the next section.

## III. Stability Criteria and Event-Triggering Condition

This section describes the main theoretical claim through the following theorem.

**Theorem 1.** $\forall\, \sigma \in (0, 1)$ *system (11) with control law (13) ensures ISS with respect to its measurement error* $e_y(t)$ *if there exist an event triggering sequences* $t_{k\in I}$ *as given below:*

$$t_0 = 0, t_{k+1} = inf\{t \in \mathbb{R} | t > t_k \wedge ||e||^2 \geq \sigma\gamma^{-1}||y_1||^2\} \tag{17}$$

*Explanation of (17)*: The equation (17) defines the event occurring condition and as well as the next event-triggering instant. According to (17), the next event ($t_{k+1}$) can occur only when $t > t_k$ and the condition as described in detail in (25), is satisfied. If $E$ is a set of entire event triggering sequences. $E_1$ is subset of E represents the sequence of $t_k$ and consecutive events. Then the infimum value of $E_1$ will be the next event occurring instant.

The parameter $\gamma$ is defined in following proof.

***Proof.*** Assuming $V(x) = x^T P x$ as an ISS Lyapunov function for (11), where $P > 0$ satisfy the following algebraic Lyapunov equation

$$A^T P + PA = -Q \tag{18}$$

where $Q$ is positive definite matrix. Then the time derivative of $V\,(x)$ along the state trajectory of (11) is

$$\dot{V}(x) = \dot{x}^T P x + x^T P \dot{x} \tag{19}$$

Now substituting (11) in (19), the $\dot{V}(x)$ is written as:

$$\begin{aligned} \dot{V}(x) &= \left(Ax + Be_y\right)^T Px + x^T P\left(Ax + Be_y\right) \\ &= x^T(A^T P + PA)x + 2x^T PBe_y \end{aligned} \tag{20}$$

To simplify (20),we used the following inequality:

For any matrix $A_1$, $A_2$, and a scalar $\epsilon > 0$, we can write the following inequality [34]:

$$2A_1A_2 \leq \epsilon A_1{}^T A_1 + \frac{1}{\epsilon} A_2{}^T A_2 \tag{21}$$

Now considering $A_1 = x$ and $A_2 = PBe_y$, following inequality is achieved :

For any scalar $\epsilon > 0$, the following inequality holds

$$2x^T PBe_y \leq \epsilon||x||^2 + \frac{1}{\epsilon}||B^T PPB||\,||e_y|| \tag{22}$$

Substituting the upper bound of $2x^T PBe_y$ in (20), the $\dot{V}(x)$ reduces to

$$\begin{aligned} \dot{V}(x) &\leq x^T(A^T P + PA)x + \epsilon x^T x + \frac{1}{\epsilon} e_y{}^T B^T PPBe_y \\ &= \lambda_{min}(Q)||x||^2 + \epsilon||x||^2 + \tfrac{1}{\epsilon}||B^T PPB||\,||e_y\,||^2 \end{aligned} \tag{23}$$

Selecting $\epsilon = \lambda_{min}(Q)/2$ , the above equation turns as

$$\dot{V}(x) \leq -\frac{\lambda_{min}(Q)}{2}||x||^2 + \frac{2}{\lambda_{min}(Q)}||B^T PPB||\,||e_y\,||^2 \tag{24}$$

In (24), the first term depends on augmented system states $x$, but at the location of generator 2, the output $y_1$ is only available. Therefore, the triggering condition should be written in terms of $y_1$, instead of $x$. To achieve that, output $y_1$ can be expressed as $y_1 = \tilde{C}x$, where matrix $\tilde{C} \in \mathbb{R}^{(n+m)\times 1}$. With this observation, the event will be triggered only when the following condition is satisfied.

$$||e_y||^2 \geq \frac{\sigma}{\gamma}||y_1||^2 \tag{25}$$

where $\sigma$ and $\gamma$ are design parameters. The value of $\gamma$ can be selected from (24), which will ensure the ISS of (11) with respect to error $e_y$ analytically. For $\gamma = \frac{[\lambda_{min}(Q)]^2}{4||B^TPPB||}$ with the event-triggering condition,

$$||e_y||^2 \geq \sigma \frac{[\lambda_{min}(Q)]^2}{4||B^TPPB||}||\tilde{C}^T\tilde{C}||\,||x||^2 \tag{26}$$

(24) reduces to

$$\dot{V}(x) \leq (\sigma||\tilde{C}^T\tilde{C}||-1)\,\lambda_{min}(Q)\,||x||^2 \tag{27}$$

Therefore $\forall\,\sigma$ which makes $(\sigma||\tilde{C}^T\tilde{C}||-1) \leq 0$, ensures the asymptotically stability of (11) as well as (12) through the control triggering law (13) and (26). In this problem the matrix $\tilde{C}$ is selected in such a way that $\left||\tilde{C}^T\tilde{C}|\right| = 1$, which basically indicates that the value of $\sigma$ may varies from (0, 1).

### A. *Inter-event time*

In event-triggered mechanism, it is essential to show that the inter-event time ($\tau = t_{k+1} - t_k > 0$) is always positive, which nullify the effects of so called Zeno behavior [35]. From event-triggering condition (25), the required time to evolve $\frac{||e_y||}{||y_1||}$ from 0 to $\gamma$ decides the inter-event time, $\tau$. To find the exact analytical expression of the $\tau$, an attempt is made in the form following theorem [14], [36].

**Theorem 2.** $\forall\,\sigma\,\epsilon\,(0,1)$, *event-triggered system (15) with event-triggered law (25), there exist an event triggering sequences* $t_{k\in I}$

$$t_0 = 0,$$

$$t_{k+1} = inf\ (t \in R|t > t_k \wedge ||e_y||^2 \geq \sigma\gamma^{-1}||y_1||^2) \tag{28}$$

*which ensures that* $\tau > 0$.

***Proof :*** From [31], the time-derivative of $\frac{||e_y||}{||y_1||}$ is

$$\frac{d}{dt}\frac{||e_y||}{||y_1||} = \frac{e_y^T\dot{e}_y}{||e_y||\,||y_1||} - \frac{y_1^T\dot{y}_1||e_y||}{||y_1||||y_1||||y_1||}$$
$$\leq \frac{||e_y||||\dot{e}_y||}{||e_y||||y_1||} + \frac{||y_1||||\dot{y}_1||||e_y||}{||y_1||||y_1||||y_1||}$$

$$= \left(1 + \frac{||e_y||}{||y_1||}\right)\frac{||\dot{e}_y||}{||y_1||} \tag{29}$$

Using triangular inequality of vector norm on (14) can be written as:

$$||\dot{e}_y|| \leq ||\tilde{C}L||||x|| + ||\tilde{C}L||||e_y|| \tag{30}$$

Where $y_1 = \tilde{C}x$. Applying (30) in (29), the inequality (29) can be simplified as:

$$\frac{d}{dt}\frac{||e_y||}{||y_1||} \leq \left(1 + \frac{||e_y||}{||y_1||}\right)\left\{||L|| + ||\tilde{C}L||\frac{||e_y||}{||y_1||}\right\} \tag{31}$$

After simplification of (31) and by considering $M_1 = ||L||$, $M_2 = ||\tilde{C}L||$ and $M_3 = M_1 + M_2$, the (31) becomes as:

$$\frac{d}{dt}\frac{||e_y||}{||y_1||} \leq M_1\left(\frac{||e_y||}{||y_1||}\right)^2 + M_2\left(\frac{||e_y||}{||y_1||}\right) + M_3 \tag{32}$$

Using comparison Lemma, the solution of (32) should satisfy the inequality $\frac{||e_y||}{||y_1||} \leq \eta(t)$ [30]. Here $\eta(t)$ is the solution of the following differential equality:

$$\dot{\eta}(t) = M_1\eta^2(t) + (M_1 + M_2)\eta(t) + M_2 \tag{33}$$

Solving (33), the analytical expression of inter-event time $\tau$ is

$$\tau = \frac{2}{\sqrt{M_3^2 - 4M_1M_2}}\ln\left\{||\frac{2M_2\gamma + M_3}{\sqrt{M_3^2 - 4M_1M_2}}||\,||\frac{M_3 + \sqrt{M_3^2 - 4M_1M_2}}{M_3 - \sqrt{M_3^2 - 4M_1M_2}}||\right\} \tag{34}$$

From (34), the inter-event time $\tau > 0$, which ensure that there must be a finite amount of aperiodic time gap between two consecutive events.

## IV. Case Studies

### A. *Two area four machine system*

#### 1) *Test System Description*

The Kundur's four machine two area system, Fig.1, is taken as one of the test system for demonstration of the event triggering mechanism in wide area controller. It consists of two symmetrical areas linked by 230 KV two lines of 220 km length and rest of parameters can be seen in [1],[6]. For stability assessment in this system, a self-clearing fault of *0.133 sec* duration on bus 8 is used for time-domain simulations. The power transfer from area-1 to area-2 is nearly 500 MW and one of the tie line circuit is removed, to make system more stressed. All generators are equipped with their local PSS using the rotor speed as its input signal. The transfer function of the PSS is shown in (35) and its output is limited by ±0.15 pu [33],[37].

$$TF_{LPSS}(s) = 30\frac{10s}{10s+1}\left(\frac{0.05s+1}{0.03s+1}\right)\left(\frac{3s+1}{5.4s+1}\right) \tag{35}$$

An additional lead-lag(L-L) WADC was designed as mentioned in section II.

#### 2) *Signal selection and control location of WADC*

The residue method has been used to identify the most appropriate controller location as well as the optimal remote signal [23],[33]. The modal residues for the exciter input and the machine speed output have been calculated. The machine 3 (in area 2) is found to be the best location for WADC with its input as the difference in the speed of machine 3 and 2. Speed deviation of machine 2 is the remote input signal to which event triggering mechanism is applied.

### B. *IEEE 39 bus 10 machine power system*

#### 1) *Test System Description*

A larger multi-machine system was also used to demonstrate the effectiveness of the event triggered WADC. This system consists of 39 buses and 10 machines as shown in Fig. 2. All synchronous machines are equipped with static excitation system, governor except machine 1 and 2. The machine 1 is an equivalent unit of large inertia and machine 2 is behaving like load [38]. The local PSS were placed as per [28]. The local lead-lag type PSS were tuned using Particle swarm Optimization (PSO), by considering the Integral Square Error (ISE) of rotor angle deviations, as an objective function[39],[40].

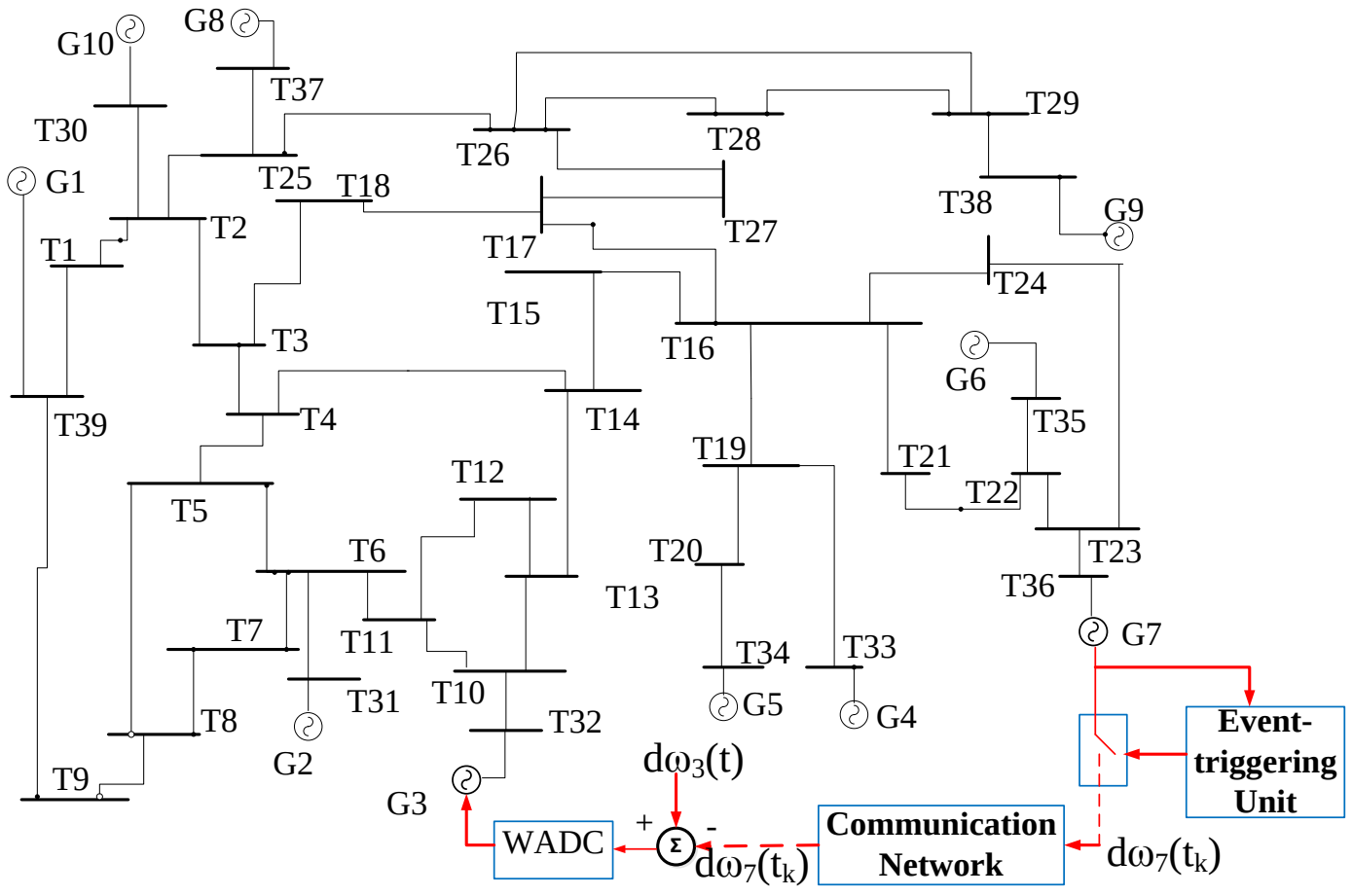

Fig. 2. *IEEE 39 bus 10 machine systemwith event-triggering mechanism*

*2) Signal selection and control location of WADC*

The WADC and the origin of its remote signal was decided on the basis of the controllability and observability index of the dominant inter area modes. The machine 3 was found to be best suited for the location of the WADC and the remote signal for the same was identified as speed deviation signal of machine 7. The complete single line diagram along-with the event triggering mechanism for WADC is shown in Fig. 2.

## V. Results and Discussions

### *A. Two area four machine power system*

The original two area test power system has four generators, and the total order of the system including the voltage regulator, turbine, governor, and local PSS for each generator, is 92. To realize the control law and for mathematical analysis, the original system is reduced to a $12^{th}$ order system using model order reduction technique available in MATLAB toolbox. The frequency response of the actual system and reduced order system is compared using Bode plot, and they are found to fairly match in the range of 1 - 100 *rad/s*, as shown in Fig. 3.

For different $\sigma$'s, the time-response of ET-WADC is compared with the base case without having WADC, as shown in Fig. 4. At this operating point, the base case system becomes unstable whenever there is a fault. However, under the same scanario if the system is equipped ET-WADC, the complete system stable. For realization of ET-WADC, the parametric value of $\sigma$ is chosen as 0.2, 0.5 and 0.9. If the value of $\sigma$ is increased beyond 0.9, the system become unstable. This demonstrates the need of WADC in the system. From Fig. 5, it is observed that increasing value of $\sigma$ and the number of event-occurrence has an inverse relation. This signifies that, parameter $\sigma$ regulates the number of transmissions. Hence, the communication burden for $\sigma$ = 0.9 is much lesser as compared to $\sigma$ = 0.2, which is also mentioned in Table I. Fig. 6 indicates that event triggering condition (25) is always satisfied during run-time and it also ensures the convergence of $||e_y||$ towards zero. The existence of positive inter-event time $\tau$, and the total number of events occurred are reported in Fig. 7. This reveals that event triggered approach drastically reduces the total number of communication over the continuous one. The same algorithmic validation has been carried out for IEEE 39-bus 10-machine system as mentioned in following subsection.

### *B. IEEE 39 bus 10 machine system*

The same procedure of order reduction described above for two area system, has been adopted for IEEE 39 bus system also. The frequency response of the actual system having order of 200 and reduced $14^{th}$ order system are compared using Bode plot, and they are found to match well in the range of inter-area oscillation frequency range. The ET-WADC for IEEE 39 bus system is implemented in the similar way as that of two area system. In IEEE 39 bus system, $d\omega_7$ is used as the remote input signal. The event-triggering law is derived using similar type of mathematical analysis as discussed in Section III. The time domain response of M1, M9 and M10 w.r.t. machine 2 is shown in Fig. 8 for different values of $\sigma$. For the sake of comparison, three values of $\sigma$ are selected as 0.1, 0.4 and 0.9. Comparison between continuous time remote signal $d\omega_7$ and event-triggered signal $d\omega_7(t_k)$ is shown in Fig. 9, which concludes the efficacy of ET-WADC over the conventional approach in terms of communication burden.

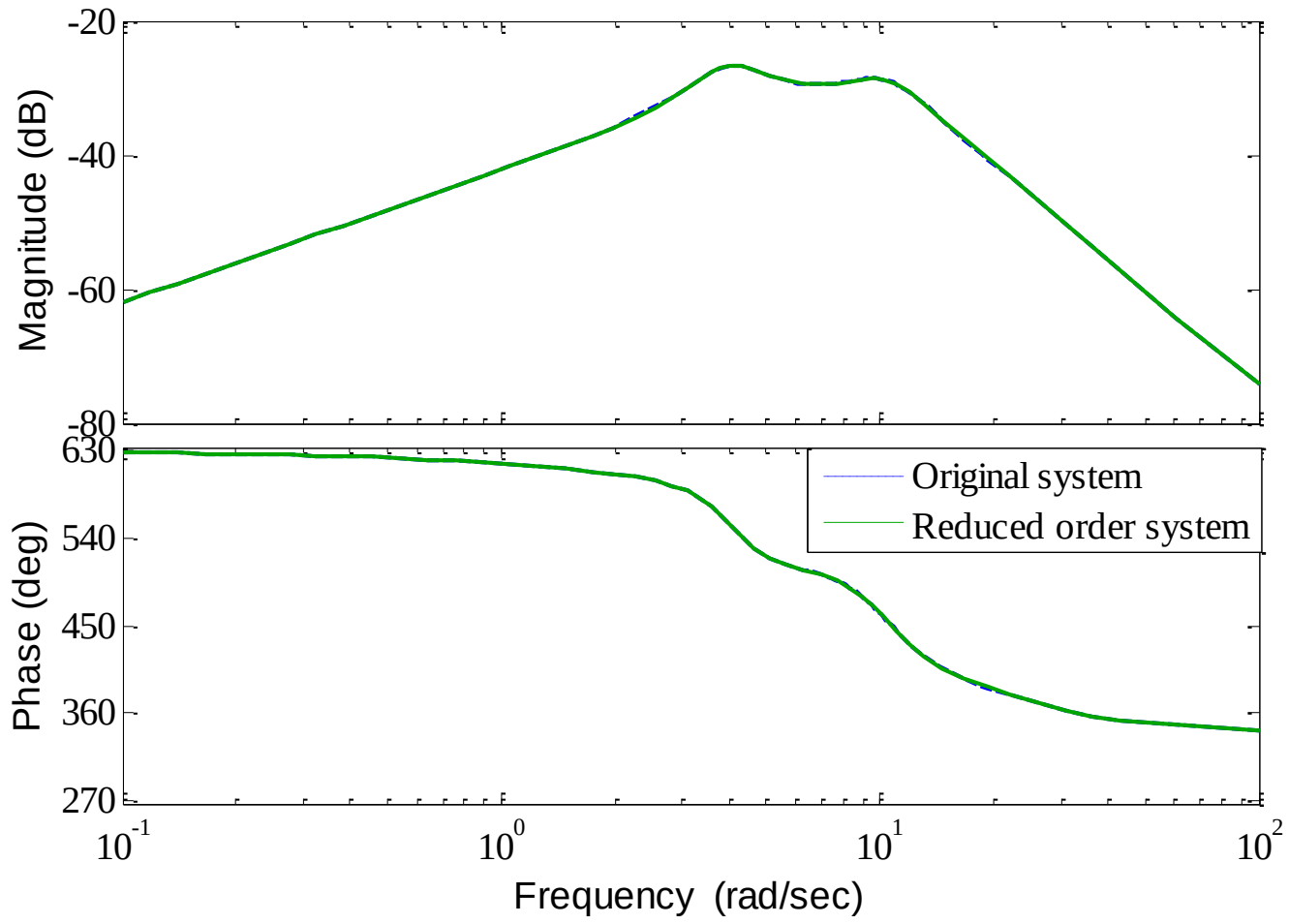

Fig. 3. Frequency response of actual and reduced order system

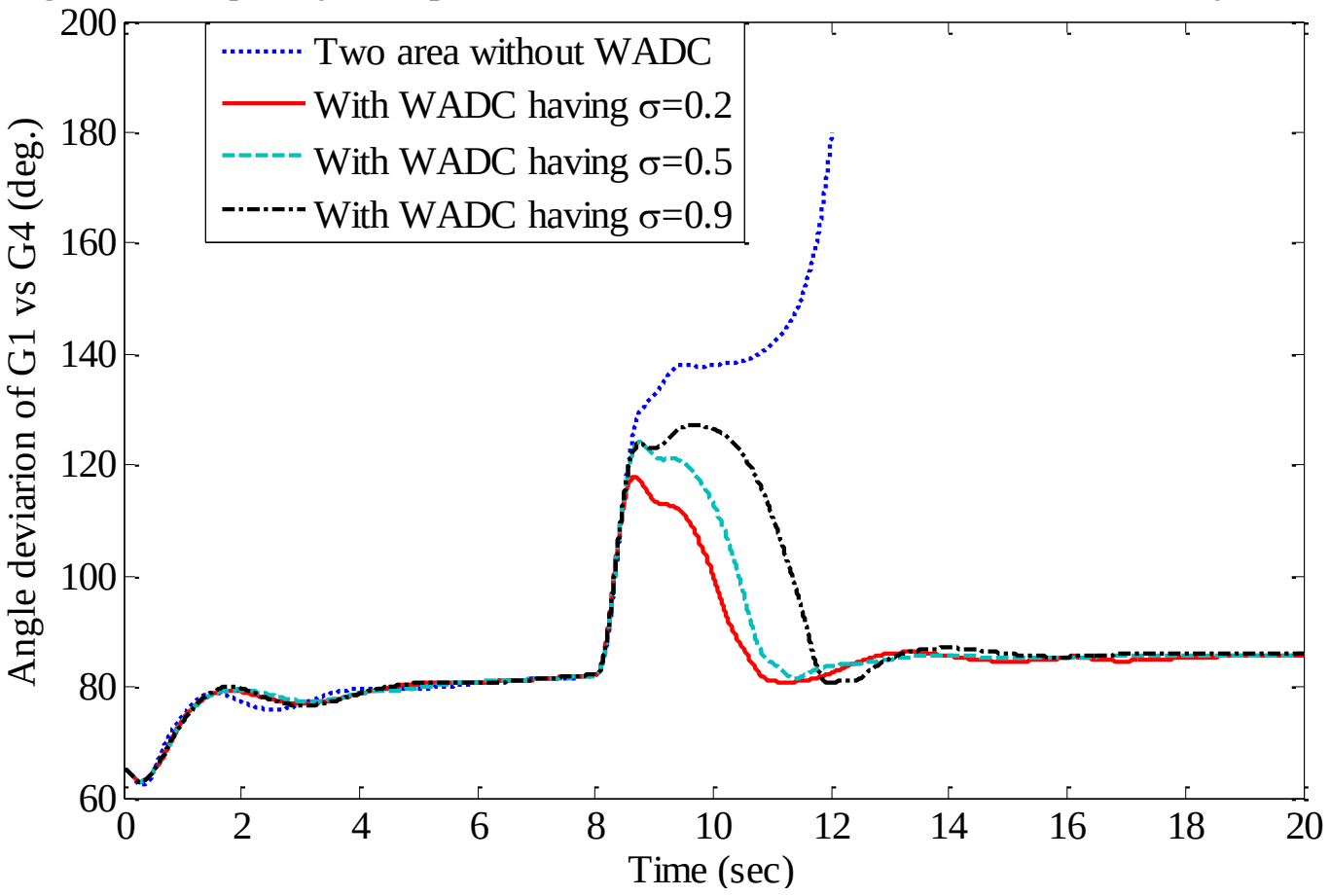

Fig. 4. *Time domain response of two-area four-machine system with and without event-triggering transmission*

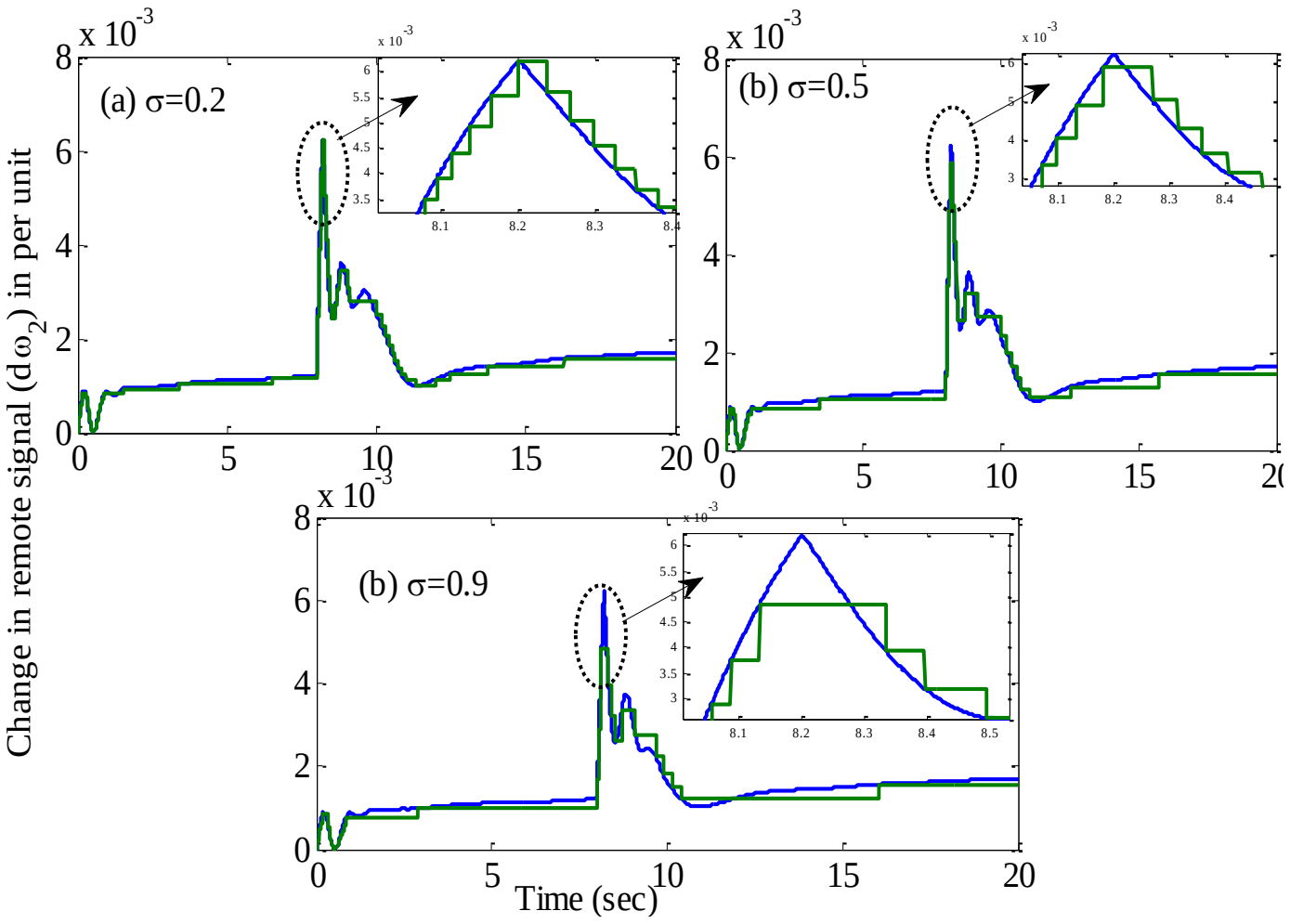

Fig. 5. *Time evolution of* $d\omega_2(t)$ *and* $d\omega_2(t_k)$ *for different values of* σ

TABLE 1:
Comparative results of ET-WADC with continuous time WADC for different σ (Two area four machine power system)

| σ | 0.2 | 0.5 | 0.9 |
|---|---|---|---|
| Continuous-time WADC | 2000 | 2000 | 2000 |
| ET-WADC | 329 | 263 | 206 |

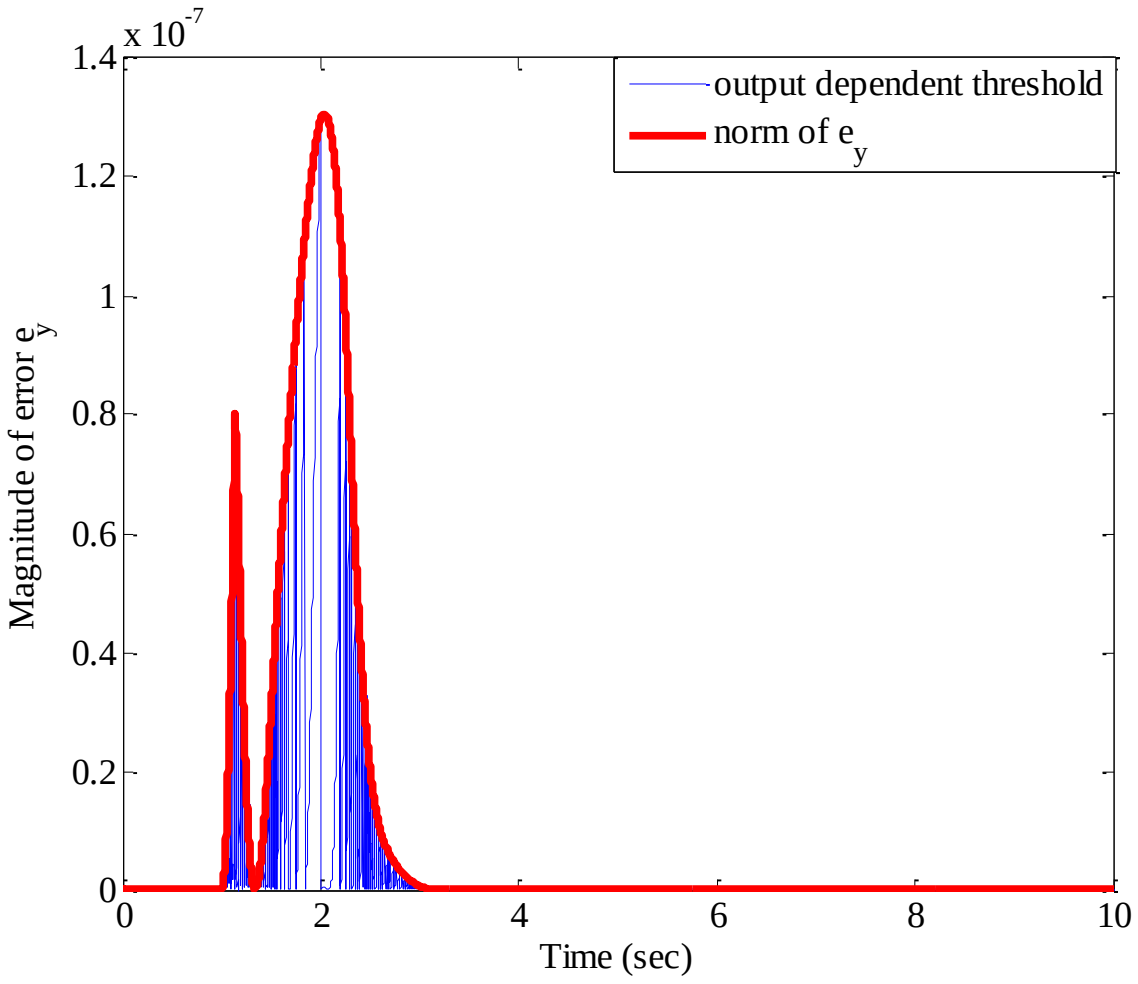

Fig. 6. The evolution of $||e_y||$ and output dependent threshold with respect to time

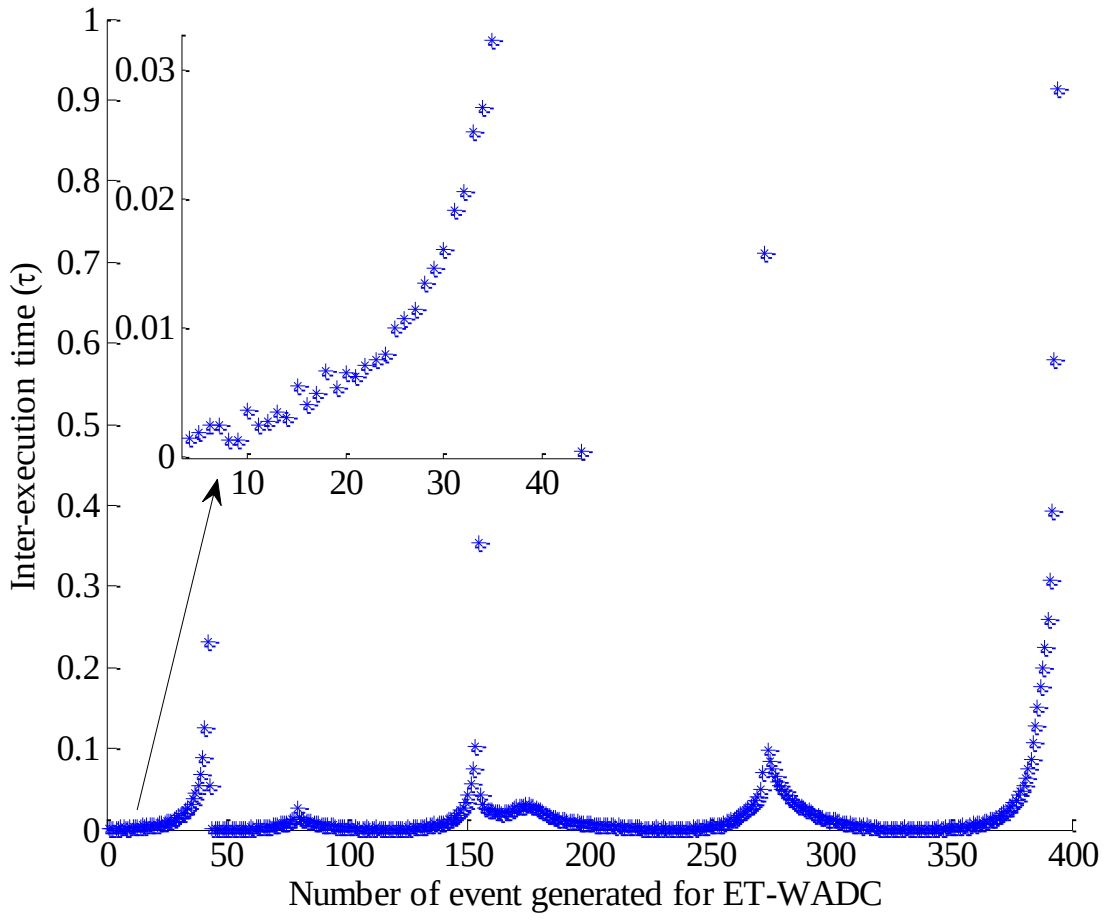

Fig. 7. Number of event-generated vs inter-event time in two-area four machine system. The y-axis represents the time-interval between two consecutive events. The x-axis denotes the total event generated during the runtime of 10 sec. The aperiodicity and positiveness of inter-event time is shown through a zoomed view (for first few seconds).

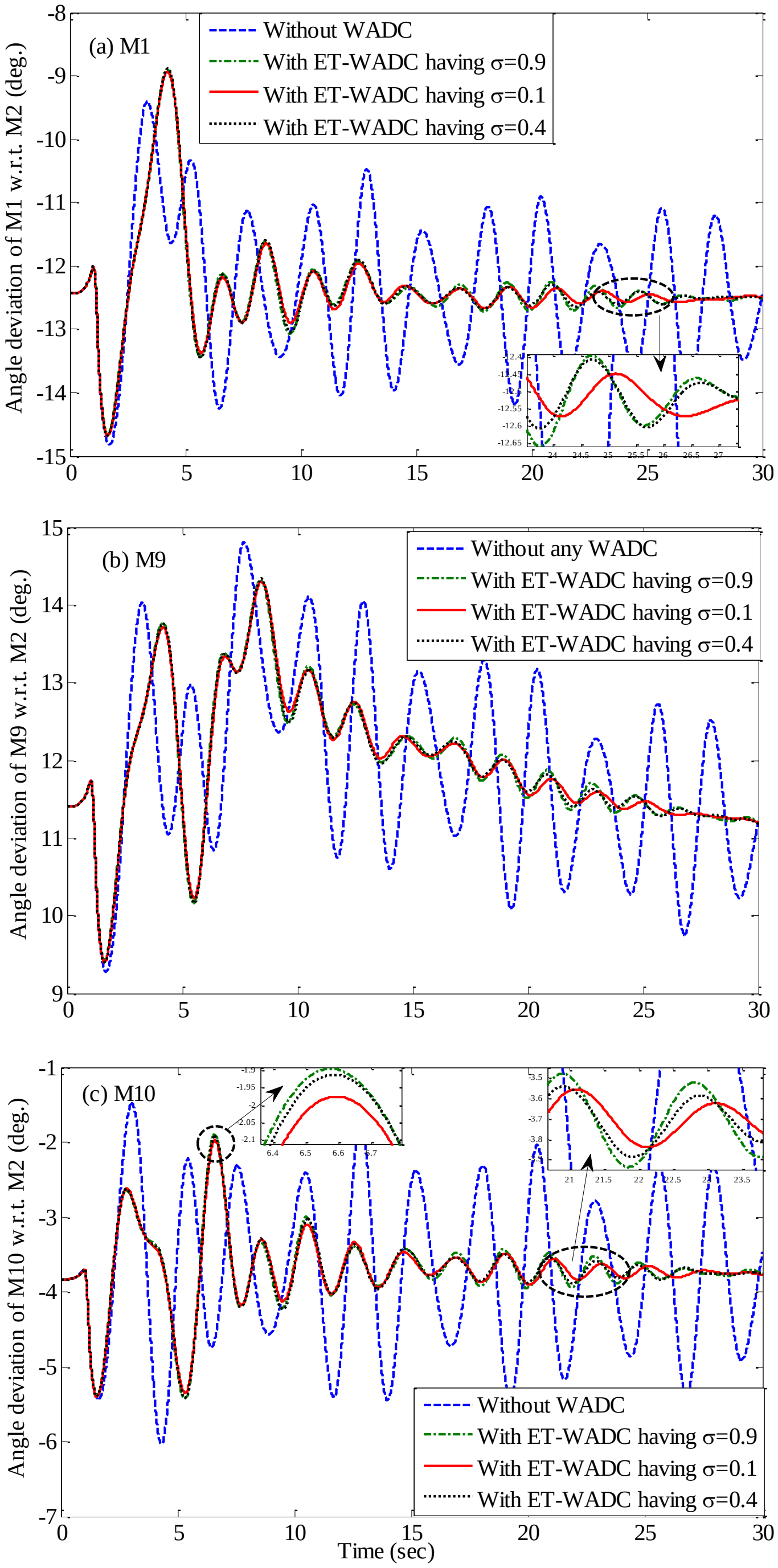

Fig. 8. Time domain response of IEEE 39 bus 10 machine system with and without event triggering transmission. (a) M1, (b) M9 and (c) M10.

TABLE 2:
Comparative results of ET-WADC with continuous time WADC for different σ (IEEE 39 bus 10 machine system)

| σ | 0.1 | 0.4 | 0.9 |
|---|---|---|---|
| Continuous-time WADC | 3000 | 3000 | 3000 |
| ET-WADC | 1445 | 1321 | 1120 |

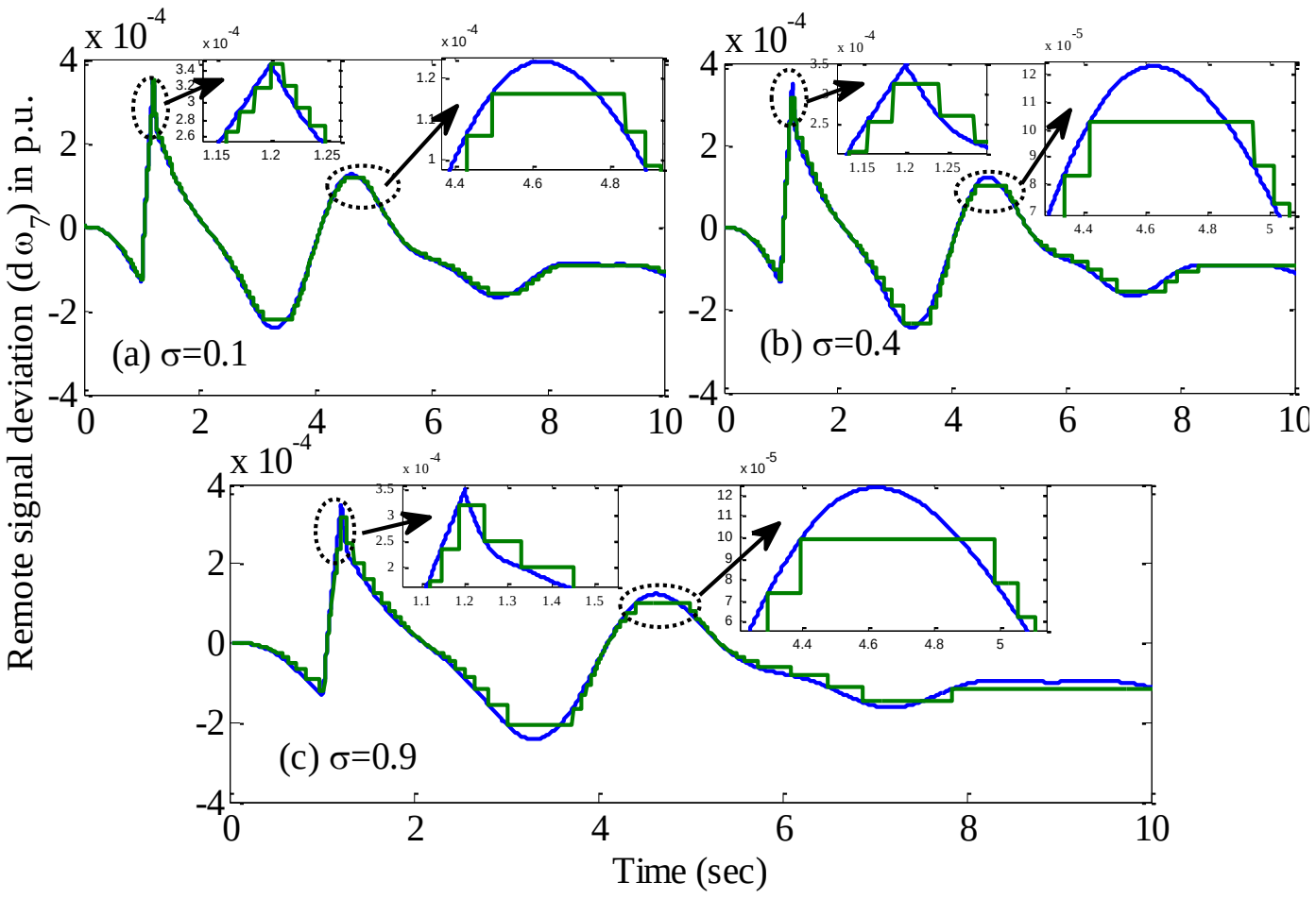

Fig. 9. *Time evolution of $d\omega_7(t)$ and $d\omega_7(t_k)$ for different values of $\sigma$*

## VI. Conclusion

This paper considers an event triggered based transmission technique for wide-area damping control system. The purpose of selecting the ET-WADC is to reduce the communication burden while maintaining damping effect of closed loop system. A trade-off between communication burden and damping effect can be established by adjusting the tuning parameter $\sigma$ in event-triggered mechanism. The asymptotic stability of closed loop event-triggered system is ensured analytically through a proper selection of ISS Lyapunov function. The effectiveness of proposed ET-WADC is demonstrated using numerical studies for two different test power systems under self clearing fault condition.

One of the primary short-coming of proposed ET-WADC is it requires continuous monitoring of event-condition. To overcome such issue, self-triggered technique can be adopted, where next event occurring instant is pre-computed based on the previous instant state or output information and measurement error [15]. Therefore the self-triggered based WADC is the future scope of this research.